\documentclass{article}
\usepackage{amssymb,amsmath,color,graphicx}
\usepackage{color}
\usepackage{hyperref}

\catcode`!=11
\let\!int\int \def\int{\displaystyle\!int}
\let\!lim\lim \def\lim{\displaystyle\!lim}
\let\!cap\cap \def\cap{\displaystyle\!cap}
\let\!cup\cup \def\cup{\displaystyle\!cup}
\let\!sup\sup \def\sup{\displaystyle\!sup}
\let\!inf\inf \def\inf{\displaystyle\!inf}
\let\!sum\sum \def\sum{\displaystyle\!sum}
\let\!max\max \def\max{\displaystyle\!max}
\let\!min\min \def\min{\displaystyle\!min}
\let\!frac\frac \def\frac{\displaystyle\!frac}
\catcode`!=12
\let\oldsection\section
\renewcommand\section{\setcounter{equation}{0}\oldsection}

\newtheorem{definition}{Definition}[section]
\newtheorem{lemma}{Lemma}[section]
\newtheorem{theorem}{Theorem}[section]
\newtheorem{remark}{Remark}[section]

\def\R{\mathbb R}
\def\Z{\mathbb Z}
\def\N{\mathbb N}
\def\epsilon{\varepsilon}
\def\ds{\displaystyle}

\newcommand\dint{\ensuremath{\ds{\int\!\!\!\!\int\hspace{-1
cm}_{_{_{_{_{_ {_{_{_{_{_{_{_{_{\R\times
(-kL,kL)}}}}}}}}}}}}}}}}}

\allowdisplaybreaks

\begin{document}

\title{ The homogenized equation of a heterogenous Reaction-Diffusion model involving pulsating traveling fronts}
\author{Mohammad El Smaily  \thanks{The author is indebted to the CMU-Portugal program and ``Center for Nonlinear Analysis'' for their support.}\\
{\footnotesize  Carnegie Mellon University,}
\\{\footnotesize Department of Mathematical Sciences, Wean Hall}
\\{\footnotesize Pittsburgh, PA, 15213, USA}\\
{\footnotesize (\href{mailto:elsmaily@andrew.cmu.edu}{elsmaily@andrew.cmu.edu})}}

\date{March 23, 2011}

\maketitle

\begin{abstract} The goal of this paper is to find the
homogenized equation of a heterogenous Fisher-KPP model in a
periodic medium.
The solutions of this model are pulsating traveling fronts
whose \emph{speeds} are superior to a parametric minimal speed
$c^*_L$.
We first find the homogenized limit of the stationary states
which depend on the space variable in many cases.
Then, we prove that the pulsating traveling fronts converge to
a classical $u_0:=u_0(t,x)$ of a homogenous reaction-diffusion
equation.
The homogenized limit $u_0$ is also a traveling front whose
minimal speed of propagation is given in terms of the
coefficients of the problem.
\end{abstract}

\bf Keywords\rm: homogenization, reaction-diffusion, front
propagation, heterogenous media.

\bf AMS Classification\rm: 35B27, 35B45, 35K55, 35K57.

\section{Introduction and Setting of the Problem}
This paper is a continuation in the study of the propagation
phenomena of pulsating traveling fronts solving a heterogenous
reaction-diffusion
equation. The notion of traveling fronts arose in 1937 in
the\textit{ homogenous} model of
Fisher \cite{fi} and Kolmogorov, Petrovsky, and Piskunov
\cite{kpp}. This model describes certain population dynamics. In
the one-dimensional case, it corresponds to the following
equation
\begin{equation}\label{kpp}
\frac{\partial u}{\partial t}=D \frac{\partial^2 u}{\partial
x^2}+ u\ (\mu-\nu
u), \ t>0, \ x\in \R.
\end{equation}
The unknown $u=u(t,x)$ is the population density at time $t$ and
position $x$, and the positive constant coefficients $D$, $\mu$
and $\nu$ respectively correspond to the diffusivity (mobility
of the individuals), the intrinsic growth rate and the
susceptibility to crowding effects.

Later, many works extended the notion of \emph{traveling
fronts} to the notion of
\textit{pulsating traveling fronts} solving a
\textit{heterogenous} reaction-advection-diffusion equation in
any dimensional space and in general periodic domains (see for
example
\cite{bh}, \cite{bhr1}, \cite{bhr2}, \cite{skt}, \cite{skt87},
\cite{w2}, \cite{x1}, \cite{x2}, \cite{x3} and \cite{x6}). We
will recall, after introducing the terms in our problem, the
definition of pulsating traveling fronts in the one-dimensional
case. The references which were mentioned above can give a
detailed and wide description of this notion in higher
dimensions and in many general settings.

In this paper, the setting is similar to that in El Smaily,
Hamel, and Roques \cite{ehr}. We consider the parametric
heterogenous reaction-diffusion equation ($L>0$ is the
parameter)
 \begin{equation}\label{eqevo}
\frac{\partial u}{\partial t}=\frac{\partial}{\partial
x}\left(a_L(x)\frac{\partial u}{\partial x}\right)+
f_L(x,u),\quad t\in\mathbb{R},~x\in\mathbb{R}.
\end{equation}
The diffusion term $a_L$ satisfies $$a_L(x)=a(x/L),$$ where $a$
is a $C^{2,\alpha}(\mathbb{R})$ (with $\alpha>0$) 1-periodic
function that satisfies
\begin{eqnarray}\label{ca2}
\exists\ 0<\alpha_{1}<\alpha_2,~\forall\
x\in\mathbb{R},~\alpha_1\leq a(x)\leq\alpha_2.
\end{eqnarray}
The reaction term satisfies
$f_L(x,\cdot)=f(x/L,\cdot)$, where we assume that
$f:=f(x,s)~:\R\times\R_+\to\R$ is 1-periodic in $x$, of class
$C^{1,\alpha}$ in $(x,s)$, $C^2$ in $s$ over $(0, M)$ ($M$ is defined below at (\ref{cf1})) and for each $x\in\R,$ $s\mapsto f(x,s)/s$ is continuous to the right at $s=0$. We set 
\begin{equation}
\mu(x):=\lim_{s\to 0^+}\frac{f(x,s)}{s} \hbox{ and }\mu_L(x):=\lim_{s\to
0^+}\frac{f_L(x,s)}{s}=\mu\left(\frac{x}{L}\right).
\end{equation}
In biological invasions, $\mu$ stands for the growth rate. Here, $\mu$ may depend on the position $x$. The more
favorable the region is, the higher the growth rate $\mu$ is. In this setting,
both $a_L$ and $f_L$ are $L$-periodic in the variable $x$.
Furthermore, we assume that:
\begin{eqnarray}\label{cf1}\left\{\begin{array}{ll}
\forall\ x\in\R,\quad f(x,0)=0, \\
\exists\ M\ge 0,\ \forall\ s\ge M,\ \forall\ x\in\R,\quad
f(x,s)\le 0.
\end{array}\right.
\end{eqnarray}
In the main result of this paper, we need the assumption
\begin{equation}\label{positivity of f}
\forall\, x\in\R,~\forall \,s\in(0,M),~f(x,s)>0.
\end{equation}
Moreover, to ensure the existence of pulsating traveling
fronts, we assume that $f$ satisfies the following condition
\begin{equation}\label{cf3}
\forall\ x\in\R,\ \ s\mapsto f(x,s)/s\hbox{ is decreasing in
}s>0.
\end{equation}
Let, for each $s\in\R,$
$$g(s):=\int_0^1f(x,s)dx\,=\,<\!f(\cdot,s)\!>_A$$
($<\!\cdot\!>_A$ stands for the arithmetic mean of a function, see Definition \ref{mean} below).
Then (\ref{cf1}) and (\ref{cf3}) yield that $g(0)=0,$
$g(s)\leq0$ for all $s\geq M,$ and
$s\mapsto \displaystyle{\frac{g(s)}{s}}$ is decreasing in $s.$

The stationary states $p(x)$ of (\ref{eqevo}) satisfy the
equation
\begin{equation}\label{eqsta}
\frac{\partial}{\partial x}\left(a_L(x)\frac{\partial
p}{\partial x}\right)+f_L(x,p)=0,~x\in\mathbb{R}.
\end{equation}
Under general hypotheses including those of this paper, and in
any space dimension, it was proved in \cite{bhr1} that a
necessary and sufficient condition for the existence of a
positive and bounded solution $p$ of (\ref{eqsta}) was the
negativity of the principal eigenvalue $\rho_{1,L}$ of the
linear operator
$$\mathcal{L}_0:  \Phi \mapsto-(a_L(x)\Phi')'-\mu_L(x)\Phi,$$
with periodicity conditions $\Phi(x)$ is $L$-periodic in $x$. In this case, the solution $p$ was
also proved to be unique, and therefore $L$-periodic. Actually,
it is easy to see that the map $L\mapsto \rho_{1,L}$ is
nonincreasing in $L>0$, and even decreasing as long as $a$ is
not constant (see the proof of Lemma 3.1 in \cite{ehr}).
Furthermore, $\rho_{1,L}\to-\displaystyle{\int_0^1}\mu(x)dx$ as
$L\to 0^+$. In this paper, the assumption (\ref{positivity of
f}) yields that $\mu(x)$ is positive everywhere and hence
\begin{equation}\label{hypl1}
\displaystyle{\int_0^1}\mu(x)dx>0.
\end{equation}
The hypothesis (\ref{hypl1}) guarantees  that (for more details about this, see (\ref{strict neg of eigen values}) below)
$$\forall\ L>0,\quad \rho_{1,L}<0.$$

\noindent Throughout this paper,  we set $<\cdot>_A$ as follows.
\begin{definition}[Arithmetic Mean]\label{mean}
Let $w:\R\mapsto\R$ be an $l-$periodic function (for some $l>0$). If $w$ is integrable over its period $[0,l]$, we define the arithmetic mean of $w$ as
\begin{equation}\label{Amean}
<w>_A:=\frac{1}{l}\int_{0}^l w(x)dx.
\end{equation}
\end{definition}

\noindent Now, we recall the definition of pulsating traveling fronts in
the one-dimensional case.

\begin{definition}[Pulsating traveling fronts]\label{pulsating
traveling fronts}
A function $u\,=\,u(t,x)$ is called a pulsating traveling front
propagating from right to left with an effective speed
$c\,\neq\,0,$ if $u$ is a classical solution of:
\begin{eqnarray}\label{front}\left\{\begin{array}{ll}
\displaystyle{\frac{\partial u}{\partial t}=
\frac{\partial}{\partial x}\left(a_L(x)\frac{\partial
u}{\partial x}\right)+f_L(x,u),\quad
t\in\mathbb{R},~x\in\mathbb{R},}\\
\displaystyle{\forall\, k\in\mathbb{Z},\;
\forall\,(t,x)\,\in\,\mathbb{R}\,\times\,\mathbb{R}},
\quad\displaystyle{u(t+\frac{kL}{c},x)\,=\,u(t,x+kL)} \hbox{,}
\\
0\,\leq\,u(t,x)\leq\,p_L(x),\\
\displaystyle{\lim_{x\,\rightarrow\,-\infty}u(t,x)\,=0\;\hbox{and}\;
\lim_{x\rightarrow \,+\infty} u(t,x)-p_L(x)\,=\,0}
\hbox{,}\end{array}\right. \label{defpuls}
\end{eqnarray}
where the above limits hold locally in $t.$
\end{definition}

This definition was given in any space dimension in \cite{bh}
and \cite{x1} whenever the stationary state $p_L\equiv1$ and in
\cite{bhr2} whenever $p_L\not\equiv 1.$

 For each $L>0$,  assuming (\ref
{ca2}) on the diffusion $a,$ (\ref{cf1}), (\ref{cf3}) and
(\ref{hypl1}) on the nonlinearity $f,$ the results of
\cite{bhr2} yield that there exists $c^{*}_L>0$ such that
pulsating traveling fronts of
(\ref{front}) propagate with a speed $c$ exist if and only if
$c\geq c^{*}_L.$ The value $c^{*}_L$ is called the \emph{minimal
speed of propagation}. We refer to
\cite{bhn1,bhn2,El Smaily,El Smaily min
max,ehr,EK1,he8,llm1,llm2,rh1,nad2,nrx,nx,Zlatos
Ryzhik,w2,zlatos} for further results on the existence and
properties of the minimal speed in the KPP case. We mention that the limit of the miniamal wave speeds was considered in \cite{llm1} but there are no results about the homogenized equation. 

In El Smaily, Hamel, Roques \cite {ehr}, the homogenized
speed was found by calculating the limit of $c^*_L$ as
$L\rightarrow 0^+$. Precisely, Theorem 2.1 in \cite {ehr} yields
that
\begin{equation}\label{limit as L tends to 0}
\displaystyle\lim_{L\rightarrow0^{+}}c_{L}^{*}=2\sqrt{<\!a\!>_{H}\
<\!\mu\!>_A},
\end{equation}
where $$<\!\mu\!>_A\ =\ \int_{0}^{1}\mu(x)dx\ \hbox{ and }\
<\!a\!>_{H}\ =\ \left(\int_{0}^{1}(a(x))^{-1}dx\right)^{-1}=\
<\!a^{-1}\!>_A^{-1}$$
denote the arithmetic mean of $\mu$ and the harmonic mean of $a$
over the interval $[0,1].$ This result was proved rigorously and
it generalized the formal and numerical results of \cite{kkts}.

Having (\ref{limit as L tends to 0}), there arise several
questions about the homogenized equation of (\ref{front}) and
the nature of the homogenous limit of the pulsating traveling
fronts $u_L$ and the type of convergence of
$\displaystyle{\left\{u_L\right\}_L}$ as the periodicity
parameter $L\rightarrow0^+.$ The main goal of this work is to
answer these questions. 

In this paper, some difficulties arise in finding
$H^1_{loc}(\R\times\R)$ estimates, independent of $L$, for a
sequence of pulsating traveling fronts
$\displaystyle{\left\{u_L\right\}_L}$ and for the corresponding
sequence $\displaystyle{\left\{a_L\,\frac{\partial\,
u_L}{\partial x}\right\}_L}$. In fact, each pulsating traveling
front $u_L$ satisfies a sort of $(t,x)$-periodicity (see the
second line of (\ref{front})). This fact makes the procedure
leading to the desired estimates indirect. Another difficulty
comes from the dependance of the stationary states $p_L$ on the
space variable $x.$ This is due to the choice of a wider class
of heterogeneous nonlinearities in the present work. We mention
that the situation becomes simpler if we assume that there is a
positive value $s_0$ such that $f(x,s_0)=0$ for any $x\in\R$ and
that $f(x,s)>0$ in $\R\times(0,s_0)$. Indeed, this yields that
$p_L\equiv s_0$ for all $L>0$ (see \cite{bhr1} and \cite{bhr2}
for more details). One of the techniques that we use in this present work appears in Step 3 of the proof of Theorem \ref{main th}. It consists of deriving the reaction-diffusion equation with
respect to the time variable and then getting estimates on the functions $w_L:=\partial u_L/\partial t$ and $v_L:=a_L(x)\partial u_L/\partial x.$ 

\section{Main Results}\label{main results}

Before going further in this section, we recall that the
function $g$ defined by
$$\forall\,s\in\R,~g(s):=\int_{0}^1f(x,s)dx$$ satisfies
$g(0)=0.$ Referring to Definition \ref{mean}, one can rewrite  
\begin{equation}
\forall\,s\in\R,~g(s)=<\cdot,s>_A.
\end{equation}
 Moreover, (\ref{hypl1}) yields that
$$g'(0)=\lim_{s\rightarrow0+}\int_{0}^{1}
\frac{f(x,s)}{s}=\int_{0}^{1}\mu(x)dx\;>0.$$ Owing to
(\ref{cf1}) and (\ref{cf3}), the map $s\mapsto\frac{g(s)}{s}$ is
decreasing and $g(s)\leq0$
for all $s\geq M$. 

\noindent \textit{As a consequence, we let $p_0$ stand for the only root of $g(s)$ which is strictly positive}. 
\vskip0.35cm

The following lemma gives several convergence results of the
sequence $\{p_L\}_{L>0}$ of stationary states as
$L\rightarrow0^+:$

\begin{lemma}[The homogenized stationary state at
$+\infty$]\label{hom pL}
Assume that the diffusion $a=a(x)$ satisfies (\ref{ca2}) and the
nonlinearity $f$ satisfies (\ref{cf1}) and (\ref{cf3}) together
with $\int_{0}^1\mu(x)dx>0$. Let $\{L_n\}_{n \in \N}$ be a
sequence of positive real numbers in $(0,1)$ such that
$L_n\rightarrow 0^+$ as $n\rightarrow+\infty.$ Let $p_0$ denote
the unique positive zero of the function $g(s):=<f(\cdot,s)>_A$ (see Definition \ref{mean} and the explanation above about $p_0$).
For each $n\in\N,$ the function $p_{L_n}:=p_{L_n}(x)$ denotes
the unique stationary state at $+\infty$ of the equation
(\ref{eqsta}) with $L=L_n$.
Then,

i) The sequence $ \{ p_{L_n}\}_{n\in \N}$ is bounded in
$H^1_{loc}(\R)$.\vskip 0.25 cm

ii) $p_{L_n}\rightharpoonup p_0$ in $H^1_{loc}(\R)$ weak and
$p_{L_n}\rightarrow p_0$ in $L^2_{loc}(\R)$ strong as
$n\rightarrow+\infty$.\vskip 0.25 cm

iii) $p_{L_n}\rightarrow p_0$ in $C^{0,\delta}_{loc}(\R)$ as
$n\rightarrow +\infty$ for all $0\leq\delta<1/2$.
\end{lemma}
\begin{remark}
We mention that the assumption (\ref{positivity of f}) is not
needed in Lemma \ref{hom pL}. Actually, concerning the
nonlinearity $f$, we assume that $\int_{0}^1\mu(x)dx>0$ in order
to guarantee the existence and the uniqueness of the stationary
state $p_L$ solving (\ref{eqsta}) for each $L>0.$ The results of
Lemma \ref{hom pL} hold in the cases where the sign of $\mu$ may
be positive in some regions (favorable regions) and negative in
others (unfavorable regions) provided that
$\int_{0}^1\mu(x)dx>0$.
\end{remark}
Now, we announce the homogenized equation of (\ref{front}) and
some convergence results of $\{u_L\}_{L>0}$ as $L\rightarrow0^+$
in the following theorem.
\begin{theorem}[Homogenized equation after normalization]\label{main th}
Assume that the diffusion $a$ satisfies (\ref{ca2}) and the
reaction $f$ satisfies (\ref{cf1}-\ref{positivity of f}) and
(\ref{cf3}). Let $\{L_n\}_{n\in \N}$ be a sequence of positive
numbers in $(0,1)$ such that $L_n\rightarrow0^+$ as
$n\rightarrow +\infty.$ For each $n\in\N,$ let
$(c_{L_n},u_{L_n})$ be \underline{the} pulsating traveling front solving
(\ref{front}) for $L=L_n$, propagating with the speed
$c_{L_n}\geq c^*_{L_n},$ and satisfying the normalization  
\begin{equation}\label{startnorm}
\forall 0<L_n\leq L_0,~~\ds{\int\!\!\!\int_{(0,1)\times(0,1)}u_{L_n}(t,x)\ \!dt\ \!dx=\frac{p_0}{2}}
\end{equation}
(a justification of this normalization is given in Step 2. of the proof).
Assume that $\{c_{L_n}\}_{n\in\N}$ converges and call
$\ds{c:=\lim_{n\rightarrow+\infty}\;c_{L_n}}$.  On the other hand, let $u_0(t,x)=U_0(x+ct)$ denote the traveling front propagating from
right to left with the speed $c,$ normalized by 
\begin{equation}\label{u0normalized}
\int\!\!\!\int_{(0,1)^2}u_0(t,x)\ \!dt\ \!dx=\frac{p_0}{2},
\end{equation}
 and which is a classical solution of the homogenous reaction-diffusion equation
\begin{equation}\label{homogenized eq}
\frac{\partial u_0}{\partial t}=<a>_H
\frac{\partial^{2}u_0}{\partial x ^2}+g(u_0)\hbox{ {\rm in }
}\R\times\R,
\end{equation}
with $U_0(-\infty)=0$ and $U_0(+\infty)=p_0$ in $C^2_{loc}(\R)$.Then,
\begin{equation}\label{mainresult}
u_{L_n}\rightarrow u_0 \text{ as } n\rightarrow+\infty \text{ in }H^{1}_{loc}(\R\times\R) \text{ weak and in }L^2_{loc}(\R\times\R)\text{ strong.}
\end{equation}
\end{theorem}

The above theorem holds for pulsating traveling fronts with speeds $c_L\in[c^*_L,+\infty).$ Physically and biologically, the fronts propagating with a speed $c^*_L$ (the minimal speed) are the most interesting. Due to 
\cite{ehr}, any sequence $\left\{c^*_{L_n}\right\}_n$ will be convergent to $2\sqrt{<a>_H}\sqrt{<\mu>_A}$ as $L_n\rightarrow 0^+.$ However, for the sake of completeness, we announced Theorem \ref{main th} for any sequence $c_L\geq c^*_L$ provided
that the chosen sequence of speeds converges to $c.$ Here we recall the existence results of pulsating traveling fronts in \cite{bhr2} together with the uniqueness results, up to a shift in the time variable, proved by Hamel and Roques \cite{hr2}. In this present setting, this reads as: for each $c_L\geq c^*_L,$ there exists a pulsating traveling front $u_L(t,x):=\varphi_L(x+c_Lt,x)$ and any other pulsating traveling front with the same speed should have the form $u_L(t+\sigma,x)$ for some shift $\sigma\in\R.$ Of course, one has to take the dependence of the shifts on
$L$ into consideration. In this context, we have the following Remark \ref{rem}.  
\begin{remark}\label{rem}
In Theorem \ref{main th}, the normalization (\ref{startnorm}) of the sequence $\{u_{L_n}\}_n$ is to guarantee the uniqueness of the limit $u_0$ in the effective equation (\ref{homogenized eq}). The estimates done in Step 3 of the theorem's proof do not require the normalization. For more details, we refer the reader to the proof of the theorem and especially Step 2 which deals with the normalization issue and Step 4 which is the passing to the limit. On the other hand, if one picks a sequence of  speeds $c_L\rightarrow c> \lim_{L\rightarrow 0^+}{c^*_L}$, Step 2 of the proof, yields that there will exist a \underline{normalized} sequence 
$\{u_{L_n}\}_n$ satisfying  (\ref{startnorm}) and then the convergence result (\ref{mainresult}) will still hold.
\end{remark}
The homogenization problem would be very interesting in multi-dimensions. Estimates   for the multi-dimensional problem in  cylindrical type domains  could be obtained in a similar approach to the one we use in the proof of Theorem \ref{main th} but they would not be enough. One needs get further estimates on higher order derivatives of the sequence $u_L.$ Furthermore, the  homogenized speed of such problems was proved rigorously only in the one-dimensional case  in \cite{ehr}.  A first step to obtain the homogenized  equation in higher dimensions is to homogenize the minimal speeds. This is the subject of a forthcoming paper for the author in collaboration with S. Kirsch. 

 We lastly mention that other homogenization results were
found by Caffarelli, Lee and Mellet \cite{clm1, clm2} in the case of
combustion-type nonlinearities.

\begin{remark}From the above theorem, we can recover the sharp
lower bound of  $\ds{\liminf_{L\rightarrow0^+}c^*_{L}}$ which was
proved in \cite{ehr}. That sharp lower
bound is given by $2\sqrt{<\!a\!>_{H}\ <\!\mu\!>_A}$ which is
the minimal speed of the the homogenized equation
(\ref{homogenized eq}).
\end{remark}

\section{Proofs of the announced results}\label{proofs}

\subsection*{Proof of Lemma \ref{hom pL}} The proof of Lemma
\ref{hom pL} will be divided into three steps:

\underline{\textbf{Step 1: Convergence to a constant limit
$p_*$.}}
Under the assumptions of Lemma \ref{hom pL} on $f$, it follows
from
\cite{bhr1} that  for each $L>0,$ the function
$p_L$ solving the equation
\begin{equation}\label{eq pL}
\left(a_L(x)p_L'\right)'+f_L(x,p_L)=0,~x\in\mathbb{R}
\end{equation}
is unique, positive, $L$-periodic
and $$\forall\,L>0,\;\forall\,x\in\R,\;0<p_L(x)\leq M$$
where $M$ is the constant appearing in (\ref{cf1}).

One can directly conclude from above that the sequence
$\{p_{L_n}\}_{n\in\N}$ is bounded in $L^2_{loc}(\R).$ Now, we
fix $L>0,$ multiply the equation (\ref{eq pL}) by $p_L$ and then
integrate by parts over any interval of the form $[-kL,kL]$
where $k\in \N.$ Owing to the $L$-periodicity of $p_L,$ we get
\begin{equation}\label{mult by pL and int}
\forall \,L>0,\;\;\forall\,k\in \N, \quad
\int_{-kL}^{kL}a(\frac{x}{L})\left(p_L'\right)^2dx=L\int_{-k}^{k}f(x,p_L(Lx))
dx.
\end{equation}
Consider the values of $L$ included in the interval $(0,1)$ and
let $\mathcal{K}$ be any compact interval of $\R.$ For each
$L>0,$ we denote
\begin{equation}\label{k_L}
k_L=\left[\frac{|\mathcal{K}|}{2L}\right]+1\in\N,
\end{equation} where
$|\mathcal{K}|$ stands for the Lebesgue measure of the interval
$\mathcal{K}$
and $[\cdot]$ stands for the integer part of a real number. One
consequently has $|\mathcal{K}|\leq2k_LL\leq|\mathcal{K}|+2L$
and
$\mathcal{K}\subseteq[-k_LL+mL,k_LL+mL]$ for some integer
$m\in\Z$ depending on $\mathcal{K}$ and $L$.

Owing to the $L-$periodicity of $f_L$, $a_L$ and $p_L$ with
respect to $x$ together with the assumption (\ref{ca2}) on the
diffusion $a,$ and using (\ref{mult by pL and int}), we obtain
$$\forall\,k\in \N, \;
\int_{-kL+mL}^{kL+mL}a(\frac{x}{L})\left(p_L'\right)^2dx=\int_{-kL}^{kL}a(\frac{x}{L})\left(p_L'\right)^2dx=L\int_{-k}^{k}f(x,p_L(Lx))
dx.$$
Consequently, for any compact interval $\mathcal{K}$ in $\R$, we
have
\begin{equation}\label{L2 estimate of p'L}
\forall\,0<L\leq 1, \quad
\int_{\mathcal{K}}\left(p_L'\right)^2dx\leq C(\mathcal{K}),
\end{equation}
where
$\displaystyle{C(\mathcal{K}):=\frac{|\mathcal{K}|+2}{\alpha_1}\max_{(x,s)\in[0,1]\times[0,M]}\left|f(x,s)\right|}$
is a positive constant independent of $L$ and depending on the
size of the compact $\mathcal{K}$. In other words, the sequence
$\{p_{L_n}\}_{n\in\N}$ is bounded in $H^1(\mathcal{K})$ for any
compact $\mathcal{K}\subset\R$ and this completes the proof of
part
i) of the Lemma.

Furthermore, we can conclude that there exists $p_* \in
H^1_{loc}(\R)$ such that, up to extraction of a subsequence,
$p_{L_n}\rightharpoonup p_*$ in $H^1_{loc}(\R)$ weak and
$p_{L_n}\rightarrow p_*$ in $L^2_{loc}(\R)$ strong as
$n\rightarrow+\infty$.\\
Using Sobolev injections, we have $H^1(\mathcal{K})$ is embedded
in $C^{0,1/2}(\mathcal{K}).$ Thus, the sequence
$\{p_{L_n}\}_{n\in\N}$ is bounded in $C^{0,1/2}(\mathcal{K})$.
Compact embeddings (Schauder's estimates) yield that, one can extract a subsequence, say $\{p_{L_n}\}_n$ without loss of generality, that converges in $C^{0,\delta}_{loc}(\R)$
for each $0\leq \delta<1/2$.   As a subsequence $\{p_{L_n}\}_n$ should also converge to $p_*$ in $H^1_{loc}(\R)$ weak and
$p_{L_n}\rightarrow p_*$ in $L^2_{loc}(\R)$ strong as
$n\rightarrow+\infty$. Consequently, for each $0\leq\delta<1/2,$  there is a subsequence $p_{L_n}\rightarrow
p_{*}$, as $n\rightarrow+\infty$, in $C^{0,\delta}_{loc}(\R).$ But since each function $p_{L_n}$ is
$L_n$-periodic (with $L_n\to 0^+$ as $n\rightarrow+\infty$), it
follows from Arzela-Ascoli theorem that $p_*$ has to be constant
over $\R$.
\vskip 0.35 cm
\underline{\textbf{Step 2: The constant limit $p_*$ is 
positive.}} To achieve this goal, we will compare the stationary
states $p_L$ with the principal eigenfunctions $\Phi_L$ of the
eigenvalue problem
\begin{equation}\label{linearized}
\mathcal{L}_0\Phi_L:=-(a_L(x)\Phi_L')'-\mu_L(x)\Phi_L=\rho_{1,L}\Phi_L\hbox{
in }\R,
\end{equation}
 which are $L$-periodic and positive in $\R$.

First, we divide (\ref{linearized}) by $\Phi_L$ and then we 
integrate by parts over $[0,L].$ It then follows from the
$L$-periodicity of $\Phi_L$ and the coefficients of
$\mathcal{L}_0$ that
$$\forall\, L>0,
\quad-\frac{1}{L}\int_{0}^La_L\left(\frac{\Phi_L'}{\Phi_L}\right)^2-\int_{0}^1\mu(x)dx=\rho_{1,L}.$$
 Hence,
\begin{equation}\label{strict neg of eigen values}
\forall\, L>0, \quad \rho_{1,L}\leq
\rho_1:=-\int_{0}^{1}\mu(x)dx <0.
\end{equation}
Next, due to the uniqueness up to multiplication by a nonzero
constant of $\Phi_L$, we can assume that $||\Phi_L||_{\infty}=1$
for every $L\in\R$.
Since the function $f(x,s)$ is 1-periodic in $x$ and of class
$C^1$ on $\R\times \R^+,$ one can then find $\epsilon_0>0$ such
that
\begin{equation}\label{lower linear bound}
\forall\, 0\leq s\leq \epsilon_0,\forall\, x\in\R,~
f(x,s)-\mu(x)s\geq \frac{\rho_1}{2}s.
\end{equation}
Having $0<\epsilon_0\,\Phi_L\leq \epsilon_0,$ we get from
(\ref{strict neg of eigen values}) and (\ref{lower linear
bound}) that
 \begin{equation}\label{epsilon Phi is a subsolution}
\begin{array}{ll}
-(a_L\epsilon_0\Phi_L')'-f(\frac{x}{L},\epsilon_0\, \Phi_L)&=
\rho_{1,L}\epsilon_0\Phi_L+\epsilon_0\mu_L(x)\Phi_L-f(\frac{x}{L},\epsilon_0\Phi_L)\vspace{4
pt}\\
&\leq
\rho_1\epsilon_0\Phi_L-\frac{\rho_1}{2}\epsilon_0\Phi_L\vspace{4
pt}\\
&=\frac{\rho_1}{2}\epsilon_0\Phi_L\,<0~\hbox{in}~\R,
\end{array}
\end{equation}
for all $L>0.$

Let us now fix any $L>0$ and, for simplicity, denote
$$\psi_L:=\epsilon_0\Phi_L.$$
We recall that the functions $p_L$ and $\psi_L=\epsilon_0\Phi_L$
are both positive and $L$-periodic. Hence, we can define
$$\gamma^*:=\sup\{\gamma>0,~p_L>\gamma\psi_L\}\geq0.$$
Assume to the contrary that $\gamma^*< 1.$
From the assumption (\ref{cf3}), we have
$f(x,\gamma^*\psi_L)>\gamma^*f(x,\psi_L)$ for all $x\in \R$.
Referring to (\ref{epsilon Phi is a subsolution}), the following
inequality then holds
\begin{equation}\label{gamma psi}
-(a_L\gamma^*\psi_L')'-f(\frac{x}{L},\gamma^*\psi_L)<0 ~\hbox{
in }~\R.
\end{equation}
Set $z:=p_L-\gamma^*\psi_L.$ Then $z\geq0,$ and there exists a
sequence $x_n \in \R$ such that $z(x_n)\rightarrow 0$ as
$n\rightarrow+\infty$ (by definition of $\gamma^*$).
Owing to the periodicity of $z$, one can then assume that
$x_n\in[0,L]$. Hence, up to extraction of some subsequence,
$x_n\rightarrow\overline{x}\in[0,L].$ From continuity,
$z(\overline{x})=0$. Besides, it follows from (\ref{eq pL}) and
(\ref{gamma psi}) that there exists a continuous function
$b=b(x)$ such that the nonnegative function $z$ satisfies
$$(a_Lz')'+b(x)z<0~\hbox{in}~\R.$$
The strong maximum principle implies that $z\equiv 0$; and
hence, $p_L\equiv \gamma^{*}\psi_L$. This contradicts with
(\ref{gamma psi}). Consequently, the assumption that $\gamma^*<
1$ is false; and thus, $p_L\geq \gamma^*\psi_L\geq
\psi_L=\epsilon_0 \Phi_L$ in $\R$. One then concludes that
$$\forall
L>0,~~\max_{x\in\R}p_L(x)=\max_{x\in[0,L]}p_L(x)\geq\epsilon_0||\Phi_L||_{\infty}=\epsilon_0.$$
On the other hand, the constant limit $p_*$ to which the
$L_n$-periodic functions $p_{L_n}$ converge uniformly on every
compact of $\R$ as $n\rightarrow+\infty$ ($L_n\rightarrow0^+$)
satisfies $p_*\geq\liminf_{n\rightarrow +\infty}\max_{x\in
\R}p_{L_{n}}(x).$ Therefore, $p_*\geq\epsilon_0>0$.
\vskip 0.35 cm
\underline{\textbf{Step 3: The constant limit $p_*$ is equal to
$p_0$.}} For each $L>0$, we call $$q_L(x)=a_L(x)p_L'(x),~ x\in
\R.$$ Equation (\ref{eq pL}) can be rewritten as
\begin{equation}\label{eq qL}
\forall L>0,~ q_L'+f(\frac{x}{L},p_L)=0\hbox{ in }\R.
\end{equation}
Consider any compact interval $\mathcal{K}$ of $\R$ and, for
each $L,$ let $k_L>0$ be the integer defined at (\ref{k_L}). From equation (\ref{eq qL}) one has $$\forall
\,0<L<1,~\int_{-k_LL}^{k_LL}(q_n')^2=L\int_{-k_L}^{k_L}f^2(x,p_L(Lx))dx.$$
Also, we have $0<p_L\leq M$ for all $L>0,$ where $M$ is the
constant appearing in (\ref{cf1}).
Thus, for each compact interval $\mathcal{K}$ of $\R$, there
exists a constant
$C_{1}(\mathcal{K}):=(|\mathcal{K}|+2)\max_{[0,1]\times[0,M]}{|f^2(x,s)|}$,
which depends only on $\mathcal{K}$, such that
  \begin{equation}\label{L2 estimates for q'n}
\forall \,0<L<1,
~\int_{\mathcal{K}}(q_L')^2(x)dx\leq\,C_{1}(\mathcal{K}).
\end{equation}
Having $\{L_n\}_n$ as a sequence of positive numbers in $(0,1)$
such that $L_n\rightarrow0^+$ as $n\rightarrow +\infty$, we
write $p_n=p_{L_n}$
and $q_n=q_{L_n}.$
The assumption (\ref{ca2}) together with (\ref{L2 estimate of
p'L}) yield that $\{q_n\}_n$ is bounded in $L^2(\mathcal{K})$.
Finally, the sequence is $\{q_n\}_n$ is bounded in
$H^{1}_{loc}(\R)$. Arguing as in Step 1, we can conclude that
there exists a constant $q_0$ such that $q_n\rightharpoonup q_0$
in $H^{1}_{loc}(\R)$ weak, $q_n\rightarrow q_0$ in
$L^2_{loc}(\R)$ strong, and
$q_n\rightarrow q_0$ in $C^{0,\delta}_{loc}(\R)$ for all
$0\leq\delta<1/2.$ However, $f(\frac{x}{L_n},p_n)\rightarrow
g(p_*)$ in $L^{\infty}(\R)$ weak-$*$ as $n\to+\infty$. Passing
to the limit as $n\to+\infty$ in equation (\ref{eq qL}) (where
$L=L_n$) implies that $g(p_*)=0.$
Referring to the properties of the function $g$ which are
mentioned at the beginning of Section \ref{proofs} and owing to
the positivity of the constant $p_*,$ we conclude that
$p_*=p_0$. Eventually, this completes the proof of Lemma
\ref{hom pL}.\hfill$\Box$
\vskip 0.35 cm

\subsection*{Proof of Theorem \ref{main th}}
This proof will be done in four steps:

\subsubsection*{Step 1: Recalling the lower bound for the speeds
 which proves that $c>0$.} From the results of \cite{bhr2},
each pulsating traveling front $u_{L_n}$ exists if and only if
$c_{L_n}\geq c^{*}_{L_n}.$ Moreover, for each $L>0,$ the minimal
speed $c^{*}_L$ is positive and, from \cite{bhr2} (see also
\cite{bhn2} in the case when $p\equiv 1$), it is given by the
variational formula
\begin{equation}\label{var formula}
c^{*}_L=\min_{\lambda>0}\frac{k(\lambda,L)}{\lambda}=\frac{k(\lambda^{*}_L,L)}{\lambda^{*}_L},
\end{equation}
where $\lambda^*_L>0$ and $k(\lambda,L)$ (for each
$\lambda\in\R$ and $L>0$) denotes the principal eigenvalue of
the problem
\begin{equation}\label{elliptic equation L periodic}
\left(a_{L}\phi_{\lambda,L}'\right)'+2\lambda
a_L\phi'_{\lambda,L}+\lambda
a_{L}'\phi_{\lambda,L}+\lambda^{2}a_L\phi_{\lambda,L}+\mu_L\phi_{\lambda,L}=k(\lambda,L)\phi_{\lambda,L}~\hbox{in}~\R,
\end{equation}
with $L$-periodicity conditions. In (\ref{elliptic equation L
periodic}), $\phi_{\lambda,L}$ denotes a principal
eigenfunction, which is of class $C^{2,\alpha}(\R)$, positive,
$L$-periodic and unique up to multiplication by a positive
constant.
In Section 3 of \cite{ehr}, the author proved with Hamel and Roques that the minimal speeds $c^{*}_{L_n}$ satisfy
$$\ds{c:=\lim_{n\rightarrow+\infty}c_{L_n}\geq\liminf_{n\rightarrow+\infty}c^{*}_{L_n}\geq
2\sqrt{\alpha_1<\!\mu\!>_A}>0.}$$ This gives a sharp lower bound for the sequence $\{c^{*}_{L_n}\}.$

\subsubsection*{Step 2: Normalization of $u_{L}$ and recalling standard change of variables.} For any $L>0,$ consider a speed $c_L\geq c^*_L$ and let $W_L(t,x)$ be an arbitrarily chosen pulsating traveling front for (\ref{front}) propagating with speed $c_L.$ We know from Hamel, Roques \cite{hr2} that the family of pulsating traveling fronts with a speed $c_L$ is then of the form $\{W_{L,\sigma}\}_{\sigma\in\R}$ where 
$$W_{L,\sigma}(t,x)=W_L(t+\sigma,x),~~(t,x)\in\R^2.$$ This is the uniqueness up to a shift in the time variable of traveling fronts with the same speed $c_L$. We want to pick one front of this family and use it in the homogenization procedure. Once the appropriate normalization is picked, we will drop the index $\sigma$ and continue the proof with the notation $u_L$ (or $\varphi_L$)  for our front.
 We start by recalling the standard change of variables $$\forall\, L>0,\forall\sigma\in\R,~W_{L,\sigma}(t,x):=\varphi_{L,\sigma}(x+c_Lt,x)=\varphi_{L,\sigma}(s,x),~(t,x)\in\R\times\R.$$
The $(t,x)-$periodicity (or the definition of the speed)  
$$W_{L,\sigma}(t,x+\frac{kL}{c_L})=W_L(t,x+k), \text{ for all }k\in \Z,$$ 
came, in \cite{bhr2} for example, by the construction of
$\varphi_{L,\sigma}$ (for any shift $\sigma\in\R$) which is $L-$periodic with respect to the spatial variable $x$ and solves  the equation 
\begin{equation}\label{phieqsigma}
\partial_x(a_L\partial_x\varphi_{L,\sigma})+a_L\partial_{ss}\varphi_{L,\sigma}+\partial_x(a_L\partial_s\varphi_{L,\sigma})+\partial_s(a_L\partial_x\varphi_{L,\sigma})
-c_L\partial_s\varphi_L+f_L(x,\varphi_L)=0,
\end{equation}
for all $(s,x)\in\R\times\R$.  This is (\ref{front}) in terms of the new function $\varphi_{\sigma,L}.$ Furthermore, one has $\varphi_{L,\sigma}(s,x)\rightarrow0$ as $s\rightarrow-\infty$ and $\varphi_{L,\sigma}(s,x)\rightarrow p(x)$ as $s\rightarrow+\infty$ in
$C^2(\R)$ (in fact, this follows mainly from the standard elliptic estimates and from the periodicity of $\varphi_{L,\sigma}$ with respect to $x$). 
As a consequence, it was proved in \cite{bhr2} that
$$\forall\,L>0, ~\sigma\in\R,~\lim_{t\rightarrow -\infty}W_{L,\sigma}(t,x)=0 \hbox{ and }\lim_{t\rightarrow+\infty}W_{L,\sigma}(t,x)=p(x)\hbox{ in } C^2_{loc}(\R).$$
 Then, the $(t,x)-$periodicity of the functions $W_{L,\sigma}$ leads to the limiting conditions $\lim_{x\rightarrow-\infty}W_{L,\sigma}(t,x)=0$ and
$\lim_{x\rightarrow+\infty}W_{L,\sigma}(t,x)=p(x)$ locally in $t$.

Now, we define the  function
$$\forall\, L>0,~I_L(\sigma)=\int_{(0,1)^2}W_L(t+\sigma,x)dtdx,$$
which is continuous over $\R.$ We recall that for each $L>0,$
the function $W_L$ is increasing in the first variable (time),
hence $I_L$ is increasing in $\sigma\in\R.$ Also, $I_L$ satisfies
$\lim_{\sigma\rightarrow-\infty}I_L(\sigma)=0$ and
$$\lim_{\sigma\rightarrow+\infty}I_L(\sigma)=\int_0^1p_L(x)dx\geq\min_{x\in\R}p_L(x)>\frac{p_0}{2}>0\hbox{ for all $0<L\leq L_0$},$$
for some $L_0$ small enough which exists by the $L$-periodicity of $p_L$ and the uniform
convergence of $p_L$ to $p_0$ as $L\rightarrow0^+$. It is important here to notice that $L_0$ is determined only by the uniform convergence of the stationary states $p_L(x)$ to $p_0$ which is the unique positive root of $g$. The stationary state $p_L$ is the same, even if one considers different shifts in the time variable of the function $u_L(t,x)$ (see \cite{bhr2} for more details). 
Now, we use the continuity and monotonicity of $I_L$ with respect to $\sigma$ to conclude that 
\begin{equation}\label{normalizing}
\ds{\forall 0<L\leq L_0, \exists \textbf{ unique }\sigma_L \in\R,\int\!\!\!\int_{(0,1)\times(0,1)}W_L(t+\sigma_L,x)\
\!dt\ \!dx=\frac{p_0}{2}.} 
\end{equation}
\textbf{Conclusion.} For each $0<L\leq L_0,$ given a speed $c_L\geq c^*_L,$  one can then denote by 
\begin{equation}\label{defineuL}
u_L(t,x):=W_L(t+\sigma_L,x)
\end{equation} ($\sigma=\sigma_L$) for the \textbf{unique} pulsating traveling front solving (\ref{front}) and propagating with a speed $c_L\geq c^*_L$ that satisfies the normalization
\begin{equation}\label{normalizing1}
\ds{ \forall 0<L\leq L_0,~~\int\!\!\!\int_{(0,1)\times(0,1)}u_L(t,x)\
\!dt\ \!dx=\frac{p_0}{2}.}
\end{equation}
In the same context ($\sigma=\sigma_L$), we set the corresponding $\varphi_L$  as 
\begin{equation}\label{definephi}
u_L(t,x):=\varphi_L(x+ct,x)=\varphi(s,x),
\end{equation}
which, from (\ref{phieqsigma}), satisfies the equation
\begin{equation}\label{phi eq}
\partial_x(a_L\partial_x\varphi_{L})+a_L\partial_{ss}\varphi_{L}+\partial_x(a_L\partial_s\varphi_{L})+\partial_s(a_L\partial_x\varphi_{L})
-c_L\partial_s\varphi_L+f_L(x,\varphi_L)=0,
\end{equation}
for all $(s,x)\in\R\times\R$.
\subsubsection*{Step 3: Boundedness of $\left\{u_{L_n}\right\}_n$ and $\left\{a_{L_n}{\partial_x u_{L_n}}\right\}_n$ in $H^{1}_{loc}(\R\times\R).$}
To simplify notations, we drop the $n$ and consider the family $\{(c_L,u_L)\}_{0<L<1}$ of pulsating traveling fronts solving
(\ref{front}), satisfying the normalization (\ref{normalizing1}) and such that $0<\underline{c}\leq c_L\leq \overline{c}$
for all $0<L<1$, where $\underline{c}$ and $\overline{c}$ are
two positive constants. We mention that, for the sequence
$\{(c_{L_n},u_{L_n})\}_{n\in\N}$ which we consider in Theorem
\ref{main th}, we have
$\underline{c}=2\sqrt{\alpha_1<\!\mu\!>_A}$ (see Step 1) and
  $\overline{c}=\sup_{n\in\N}c_{L_n}.$

   Since $\varphi_L(-\infty,x)=0$  and
$\varphi_L(+\infty,x)= p(x)$ in $C^2(\R)$ (for each $L>0$), it
follows then that $\nabla_{s,x}\varphi_L(-\infty,x)=0,$
$\partial_x\varphi_L(+\infty,x)=p'_L(x),$ and
$\partial_s\varphi_L(+\infty,x)=0$ uniformly in $x\in\R.$
Integrating (\ref{phi eq}) by parts over $\R\times[-kL,kL]$
(where $L>0$ and $k\in\N$) and using the $L-$periodicity of
$\varphi_L$ with respect to $x$, we then get
\begin{equation*}\dint
f(\frac{x}{L},\varphi_L(s,x))dsdx=c_L\int_{-kL}^{kL}p_L(x)dx-\int_{-kL}^{kL}a_Lp'_L(x)dx,
\end{equation*}
or equivalently
\begin{equation}\label{int phi eq}
\dint
f(\frac{x}{L},u_L(t,x))dtdx=\int_{-kL}^{kL}p_L(x)dx-\frac{1}{c_L}\int_{-kL}^{kL}a_Lp'_L(x)dx.
\end{equation}
As done in the proof of Lemma \ref{hom pL}, we take any compact
interval $\mathcal{K}\subset \R$ and we define $k_L\in\N$ as in (\ref{k_L}). We apply (\ref{int phi eq}) for $k=k_L.$ Having
$0<u_L(t,x)<p_L\leq M$ for all $(t,x)\in\R\times\R$ and owing
to (\ref{ca2}), (\ref{positivity of f}) and (\ref{L2 estimate of
p'L}), one then gets
 \begin{equation}\label{f(x/L,u) is integrable}
0<\int\!\!\!\int_{\R\times\mathcal{K}}f(\frac{x}{L},u_L(t,x))dtdx\leq
C_2(\mathcal{K})
 \end{equation}
for all $0<L<1,$ where
$$C_2(\mathcal{K}):=M(|\mathcal{K}|+2)+\frac{\alpha_2}{\underline{c}}\sqrt{C(\mathcal{K})}\sqrt{|\mathcal{K}|+2}$$
  is a constant independent of $L.$

Multiplying (\ref{phi eq}) by $\varphi_L$ and integrating by
parts over $\R\times(-kL,kL),$ we obtain
\begin{equation}\label{mult by phiL and int}
\begin{array}{ll}
\ds{\frac{c_L}{2}\int_{-kL}^{kL}p_L^2(x)dx}&=\ds{-\dint\left[a_L\left(\frac{\partial\varphi_L}{\partial
x}\right)^2+a_L\left(\frac{\partial\varphi_L}{\partial
s}\right)^2+2a_L\frac{\partial \varphi_L}{\partial
x}\frac{\partial\varphi_L}{\partial s}\right]dsdx}\vspace{4
pt}\\
&\ds{+\int_{-kL}^{kL}a_Lp'_Lp_Ldx+\dint
f(\frac{x}{L},\varphi_L(s,x))\varphi_L}\vspace{4 pt}\\
&=\ds{-\dint a_L\left(\frac{\partial u_L}{\partial
x}\right)^2dtdx}\ds{+\int_{-kL}^{kL}a_Lp'_Lp_Ldx}\vspace{4 pt}\\
&\ds{+\dint f(\frac{x}{L},u_L(t,x))u_Ldtdx.}
\end{array}
\end{equation}
Notice that the last integral in (\ref{mult by phiL and int})
converges because of (\ref{f(x/L,u) is integrable}) and $0\leq f(x/L,u_L)u_L\leq M f(x/L,u_L)$ in $\R\times\R.$ Moreover,
$$\forall\,L>0,~\left|\int_{-kL}^{kL}a_Lp'_Lp_Ldx\right|\leq\alpha_2M(2kL)^{1/2}\left(\int_{-kL}^{kL}{p'_L}^2\right)^{1/2}.$$
Consequently,
\begin{equation}\label{L2 dx uL}
\int\!\!\!\int_{\R\times\mathcal{K}}\left(\frac{\partial
u_L}{\partial x}\right)^2dtdx\leq C_3(\mathcal{K}),
\end{equation}
where
$C_3(\mathcal{K}):=\frac{1}{\alpha_1}\left[MC_2(\mathcal{K})+\alpha_2M(2+|\mathcal{K}|)^{1/2}\sqrt{C(\mathcal{K})}\right]$
is independent of $L.$

Now, we multiply (\ref{phi eq}) by $\partial_s\varphi_L$ and we
integrate by parts over $\R\times(-kL,kL).$ We notice that, from
the $L-$periodicity with respect to $x$ of the function
$\varphi_L$ and its derivatives together with the limits of
$\partial_s\varphi_L$ and $\partial_x\varphi_L$ as
$s\rightarrow\pm\infty$, we have
$$\dint\partial_x(a_L\partial_x\varphi_L)\partial_s\varphi_L=-\frac{1}{2}\dint\partial_s(a_L(\partial_x\varphi_L)^2)=-\frac{1}{2}\int_{-kL}^{kL}
a_L p_L'^{^2}dx,$$
while
$$\dint\partial_s\varphi_L\partial_x(a_L\partial_s\varphi_L)+\partial_s\varphi_L\partial_s(a_L\partial_x\varphi_L)=0.$$

 Thus,
$$c_L\dint\left(\frac{\partial\varphi_L}{\partial
s}\right)^2=-\frac{1}{2}\int_{-kL}^{kL} a_L
p_L'^{^2}dx+\int_{-kL}^{kL}F(\frac{x}{L},p_L(x))dx,$$ where
$F(y,s)=\int_{0}^{s}f(y,\tau)d\tau.$
Hence,
\begin{equation}\label{leading to L2 est of dt uL}
\dint\left(\frac{\partial u_L}{\partial t}\right)^2dtdx\leq
c_L\int_{-kL}^{kL}F(\frac{x}{L},p_L(x))dx\leq
\overline{c}\int_{-kL}^{kL} F(\frac{x}{L},p_L(x))dx.
\end{equation}
Consequently, for all $0<L<1,$
\begin{equation}\label{L2 dt uL}
\int\!\!\!\int_{\R\times\mathcal{K}}\left(\frac{\partial
u_L}{\partial t}\right)^2dtdx\leq\ds{\int\!\!\!\!\int\hspace{-1
cm}_{_{_{_{_{_ {_{_{_{_{_{_{_{_{\R\times
(-k_LL,k_LL)}}}}}}}}}}}}}}}\left(\frac{\partial u_L}{\partial
t}\right)^2dtdx\leq C_4(\mathcal{K}),
\end{equation}
where $C_4(\mathcal{K}):=
\overline{c}\,(2+|\mathcal{K}|)\max_{(x,s)\in\R\times[0,M]}F(x,s)$
is a positive constant which is independent of $L$ and depending
only on the compact $\mathcal{K}.$

Denote
$$v_L(t,x)=a_L(x)\frac{\partial u_L}{\partial x}(t,x)\;\hbox{
and }\;w_L(t,x)=\frac{\partial u_L}{\partial t}(t,x)\; \hbox{ in
}\;\R\times\R.$$
As already underlined, it follows from \cite{bhr2} ( and
\cite{bh} in the case $p_L\equiv1$) that $w_L=\frac{\partial
u_L}{\partial t}>0$ in $\R\times\R$ for each $L>0$. We shall now
establish some estimates (independent of $L$) for the functions
$v_L$ and $w_L$, in order to pass to the limit as $L\to 0^+$.

We first recall that for each $L>0, $ $\varphi_L\rightarrow p_L(x)$ in $C^2(\R)$ as $s\rightarrow +\infty$ and $\varphi_L\rightarrow 0$ in $C^2(\R)$ as $s\rightarrow -\infty$  (this was proved by construction of pulsating traveling fronts in \cite{bhr2} for example). Now, using the relation between $u_L$ and $\varphi_L$, one concludes that
for each $L>0,$ $u_L(-\infty,x)=0$ and $u_L(+\infty,x)=p_L(x)$
locally in $x$, and $w_L(\pm\infty,x)=0$ locally in $x$.
On the other hand, (\ref{L2 dt uL}) yields that
for each compact $\mathcal{K}$ and for each $L$,
$||w_L||_{L^2(\R\times \mathcal{K})}\leq \sqrt{C(\mathcal{K})}$.
Now, we differentiate (\ref{phi eq}) with respect to $t$
(actually, from the regularity of $f$, the function $w_L$ is of
class $C^2$ with respect to $x$). There holds
$$\ds{\frac{\partial w_L}{\partial
t}}=\ds{\frac{\partial}{\partial x}\left(a_L(x)\frac{\partial
\,w_L}{\partial x}\right)+ f'_u(\frac{x}{L},u_L)w_L\ \hbox{ in
}\R\times\R.}$$
Multiply the above equation by $w_L$ and integrate by parts over
$\R\times(-kL,kL)$. From (\ref{ca2}), (\ref{leading to L2 est of
dt uL}) and the fact that $0< u_L\leq M,$ it follows that
$$\ds{\dint\left(\frac{\partial w_L}{\partial
x}\right)^2dtdx\leq \frac{2kL\eta \,c_L}{\alpha_1}}$$
where $\eta$ is the positive constant defined by
$$\ds{\eta=\max_{(x,u)\in\R\times[0,M]}|f'_u(x,u)|\,\max_{x\in\R}|F(x,M)|\geq\frac{1}{2kL}\dint
f'_u(\frac{x}{L},u_L)w_L^2\ \!dt\ \!dx}>0.$$
Then, for each compact $\mathcal{K}\subset\R,$ there exists a
constant $C'(\mathcal{K})>0$ depending only on $\mathcal{K}$
such that
\begin{equation}\label{L2 est for dx wL}
\forall\ \! 0<L<1,\;\ds{\int\!\!\!\int_{\R\times
\mathcal{K}}}\left(\frac{\partial w_L}{\partial x}\right)^2dt\
\!dx\leq C'(\mathcal{K}).
\end{equation}
We pass now to the family $\{v_L\}_L.$ Actually, $v_{L}=a_{L}\frac{\partial u_{L}}{\partial x}$ and $0<\alpha_1\leq a_{L}\leq\alpha_2 $ for each $L>0.$ Thus,
(\ref{L2 dx uL}) yields that for each compact interval $\mathcal{K}$ of $\R$ and for each $0<L<1$, $\ds{||v_{L}||_{L^2(\R\times \mathcal{K})}\leq\alpha_2\sqrt{ C_3(\mathcal{K})}.}$
Furthermore,~(\ref{front}) implies that
$$\forall\, L>0,~\ds{\frac{\partial v_{L}}{\partial
x}=\frac{\partial u_{L}}{\partial t}-f(\frac{x}{L},u_{L})}\hbox{ in }\R\times\R,$$ while $0\le f(x/L,u_{L}(t,x))\leq \kappa$ in $\R\times\R$ where $\kappa=\max_{\R\times[0,M]}f(x,u)>0$ is independent of~$L.$
Together with (\ref{L2 dt uL}), one concludes that any family $\left\{{\partial v_{L}}/{\partial x}\right\}_{0<L<1}$ is
bounded in $L^2_{loc}(\R\times\R)$ by a constant independent of $L.$ On the other hand, $$\forall L>0,~\frac{\partial v_{L}}{\partial t}=a_{L}\frac{\partial^2u_{L}}{\partial t\partial x}=a_{L}\frac{\partial w_{L} }{\partial x}\hbox{ in }\R\times\R.$$ Owing to (\ref{ca2}) and (\ref{L2 est for dx
wL}), any family $\left\{\frac{\partial v_{L}}{\partial
t}\right\}_{0<L<1}$ is bounded in $L^2_{loc}(\R\times\R).$

\subsubsection*{Step 4: Passage to the limit as
$n\rightarrow+\infty$ ($L\rightarrow 0^+$).} In this step, we
consider the sequence $\{L_n\}_{n\in\N}$ of Theorem \ref{main
th} which is in $(0,1)$ and which tends to $0$ as
$n\rightarrow+\infty.$ As a consequence of the previous step,
$\{v_{L_n}\}_{n\in\N}$ is bounded in $H^1_{loc}(\R\times\R)$.
The estimates (\ref{L2 dx uL}) and (\ref{L2 dt uL}) imply that
the sequence $\{u_{L_n}\}_{n\in\N}$ is bounded in
$H^1_{loc}(\R\times\R).$ Thus, there exist $u_0$ and $v_0$ in
$H^{1}_{loc}(\R\times \R)$ such that, up to extraction of a
subsequence, $u_{L_n}\rightarrow u_0,$ $v_{L_n}\to v_0$ strongly
in $L^2_{loc}(\R\times\R)$ and almost everywhere
in $\R\times \R$,
$$\left(\frac{\partial u_{L_n}}{\partial t},\frac{\partial
u_{L_n}}{\partial x}\right)\rightharpoonup\left(\frac{\partial
u_0}{\partial t}, \frac{\partial u_0}{\partial x}\right)\hbox{
weakly in }L^2_{loc}(\R\times \R),$$
and $$\left(\frac{\partial v_{L_n}}{\partial t},\frac{\partial
v_{L_n}}{\partial x}\right)\rightharpoonup\left(\frac{\partial
v_{0}}{\partial t},\frac{\partial v_{0}}{\partial
x}\right)\hbox{ weakly in }L^2_{loc}(\R\times\R)$$
as $n\rightarrow+\infty.$ However,
$a_{L_n}^{-1}\rightharpoonup\,<a^{-1}>_A=<a>_H^{-1}\hbox{ in
}L^{\infty}(\R)$ weak-$*$ as $n\to+\infty$. Thus,
$$\frac{\partial u_{L_n}}{\partial
x}=\frac{v_{L_n}}{a_{L_n}}\rightharpoonup\frac{v_0}{<a>_H}\hbox{
weakly in }L^2_{loc}(\R\times\R)\hbox{ as }n\to+\infty.$$
By uniqueness of the limit, one gets $v_0=<a>_H\frac{\partial
u_0}{\partial x}$. On the other hand,
$f(\frac{x}{L_n},u_{L_n})\rightarrow g(u_0)$ in $L^{\infty}(\R\times\R)$
weak-$*$ as $n\to+\infty.$ Passing to the limit as $n\to+\infty$
in the first equation of (\ref{front}) with $L=L_n$ implies that
$u_0$ is a weak solution of the equation
$$\frac{\partial u_0}{\partial t}=\frac{\partial v_0}{\partial
x}+g(u_0)=<a>_H\frac{\partial^2 u_0}{\partial x^2}+g(u_0)\hbox{
in }\mathcal{D}'(\R\times\R).$$
From parabolic regularity, the function $u_0$ is then a
classical solution of the homogenous equation
$$\frac{\partial u_0}{\partial t}=<a>_H
\frac{\partial^{2}u_0}{\partial x ^2}+\,g(u_0)\hbox{ in
}\R\times\R,$$
such that $0\le u_0\le p_0$ and $\frac{\partial u_0}{\partial
t}\geq0$ in $\R\times\R$. Lastly, it follows directly from the normalization (\ref{normalizing1}) that 
\begin{equation}\label{endingnormal}
\int\!\!\!\int_{(0,1)^2}u_0(t,x)\ \!dt\ \!dx=\frac{p_0}{2}.
\end{equation}

On the other hand, it follows from the second equation of
(\ref{front}) that
$$\forall \gamma\in
\R,~~u_0(t+\frac{\gamma}{c},x)=u_0(t,x+\gamma)\ \hbox{ in
}\R\times\R,$$
where $c=\lim_{n\to+\infty}c_{L_n}>0$. In other words,
$u_0(t,x)=U_0(x+ct)$, where $U_0$ is a classical solution of the
equation
\begin{equation}\label{hom U_0}
\ds{cU'_0=<a>_HU''_0+<\mu>_Ag(U_0)~\hbox{ in }\R},
\end{equation}
that satisfies $U_0'\ge 0$ in $\R$, $0\leq U_0(s)\leq p_0$ for
all $s\in\R.$   

\noindent Standard elliptic estimates on (\ref{hom U_0}) imply that $U_0$
converges as $s\to\pm\infty$ in $C^2_{loc}(\R)$ to two constants
$U^{\pm}_0\in[0,p_0]$ such that $g(U^{\pm}_0)=0$.

\noindent Thanks to the normalization (\ref{endingnormal}), $u_0(t,x)=U_0(x+ct)$ satisfies the normalization
\begin{equation}\label{normalized U0}
\int_0^1\left(\int_{cs}^{cs+1}U_0(\tau)d\tau\right)ds=\frac{p_0}{2}.
\end{equation}

\noindent The monotonicity of $U_0$ and the nature of the function $g$ imply that $U^-_0=0$ and
$U^+_0=p_0.$ Together with the normalization (\ref{normalized U0}), one concludes that $U_0$ is the unique traveling front
for the homogenized equation (\ref{hom U_0}) with a speed $c$
and limiting conditions $0$ and $p_0$ at infinity.  We mention that, since the minimal speed for the problem (\ref{hom U_0}) is equal to
$\ds{2\sqrt{<a>_Hg'(0)}=2\sqrt{<a>_H<\mu>_A}}$, one can review that $c\geq \ds{2\sqrt{<a>_H<\mu>_A}},$ which was proved by other tools in \cite{ehr}. Eventually, the proof of Theorem \ref{main th} is complete.\hfill$\Box$

\section*{Acknowledgments}
The author would like to thank the Referees for their valuable suggestions and comments which helped  in improving the presentation of  Theorem \ref{main th} in this work.

\end{document}